\documentclass{amsart}
\usepackage{amssymb,amsthm,amsmath}

\theoremstyle{plain}

\newtheorem*{thmS}{Theorem S}
\newtheorem*{corS}{Corollary S}
\newtheorem*{thmO}{Theorem O}
\newtheorem*{thmA}{Theorem A}
\newtheorem*{corA}{Corollary A}
\newtheorem*{corA'}{Corollary A$'$}
\newtheorem*{thmB}{Theorem B}
\newtheorem*{corB}{Corollary B}

\theoremstyle{definition}
\newtheorem{remark}{Remark}

\newcommand\N{{\mathbb{N}}}
\newcommand\class{{\mathcal{C}}}
\newcommand\cprime{$'$}

\begin{document}

\title[Embeddings into simple groups]{A pastiche on embeddings into simple
groups\\ (following P.~E.~Schupp)}

\author{Zoran {\v S}uni{\'c}}

\address{Department of Mathematics, Texas A\&M University, College
Station, TX 77843-3368, USA}

\email{sunic@math.tamu.edu}

\thanks{Partially supported by NSF grant DMS-0600975}



\begin{abstract}
Let $\lambda$ be an infinite cardinal number and let $\class= \{ H_i
\mid i \in I\}$ be a family of nontrivial groups. Assume that $2
\leq |I|\leq \lambda$, $|H_i|\leq \lambda$, for $i \in I$, and at
least one member of $\class$ achieves the cardinality $\lambda$.

We show that there exists a simple group $S$ of cardinality
$\lambda$ that contains an isomorphic copy of each member of
$\class$ and, for all $H_i$, $H_{i'}$ in $\class$ with
$|H_{i'}|=\lambda$, is generated by the copies of $H_i$ and $H_{i'}$
in $S$.

This generalizes a result of Paul E.~Schupp (moreover, our proof
follows the same approach based on small cancelation). In the
countable case, we partially recover a much deeper embedding result
of Alexander Yu.~Ol{\cprime}shanski\u{\i}.
\end{abstract}

\maketitle

\section{Background and results}

In~\cite{schupp:simple} Schupp used small cancelation theory
(construction of Adian-Rabin type) to prove, among other things, the
following result.

\begin{thmS}[Schupp~\cite{schupp:simple}]
Let $G$, $H$ and $K$ be nontrivial groups with $|G| \leq |H*K|$ and
$|K| \geq 3$. There exists a simple group $S$ that contains an
isomorphic copy of $G$ and is generated by isomorphic copies of $H$
and $K$.
\end{thmS}

\begin{corS}
Let $G$ be a countable group. For all $p,q \in \{2,3,\dots\} \cup
\{\infty\}$ with $q \geq 3$, there exists a simple group $S$ that
contains an isomorphic copy of $G$ and is generated by a pair of
elements of order $p$ and $q$, respectively.
\end{corS}

The simple group constructed by Schupp in Theorem~S, in addition to
being dependent on $G$, depends on $H$ and $K$. Accordingly, the
simple group in Corollary~S, in addition to being dependent on $G$,
depends on the pair $(p,q)$.

We will show that the argument used by Schupp can be adapted in such
a way that the same simple group can be used even if one
considerably varies $H$ and $K$ in Theorem~S and, consequently, the
same simple group can be used independently of the pair $(p,q)$ in
Corollary~S.

\begin{thmA}
Let $|I| \geq 2$ and $\class= \{ H_i \mid i \in I\}$ be a countable
family of countable nontrivial groups, at least one of which has at
least 3 elements (the groups may be isomorphic for different values
of the index).

There exists a 2-generated simple group $S$ that contains an
isomorphic copy of each member of $\class$ and, for all $H_i$, $H_j$
in $\class$ with $|H_j|\geq 3$, is generated by the copies of $H_i$
and $H_j$ in $S$.
\end{thmA}

\begin{corA}
Let $G$ be a countable group.

There exists a simple group $S$ that contains an isomorphic copy of
$G$ and, for all $p,q \in \{2,3,\dots\} \cup \{\infty\}$ with $q
\geq 3$, is generated by a pair of elements of order $p$ and $q$,
respectively.

Moreover, if $|G|\geq 3$, then, for every $p \in \{2,3,\dots\} \cup
\{\infty\}$, the simple group $S$ is generated by $G$ and an element
of order $p$.
\end{corA}

In this note, countable means finite or countably infinite. The
countability limitations imposed in Theorem~A are natural since
every countable group contains only countably many finitely
generated subgroups. An extension of Theorem~S in which countability
assumptions are not used follows.

\begin{thmB}
Let $\lambda$ be an infinite cardinal number and let $\class= \{ H_i
\mid i \in I\}$ be a family of nontrivial groups. Assume that $2
\leq |I|\leq \lambda$, $|H_i|\leq \lambda$, for $i \in I$, and at
least one member of $\class$ achieves the cardinality $\lambda$.

There exists a simple group $S$ of cardinality $\lambda$ that
contains an isomorphic copy of each member of $\class$ and, for all
$H_i$, $H_{i'}$ in $\class$ with $|H_{i'}|=\lambda$, is generated by
the copies of $H_i$ and $H_{i'}$ in $S$.
\end{thmB}

\begin{corB}
For any group $G$ with $|G| \geq 3$, there exists a simple group $S$
that contains an isomorphic copy of $G$ and, for every $p \in
\{2,3,\dots\} \cup \{\infty\}$, is generated by $G$ and a single
element of order $p$.
\end{corB}

In the countable case, the embedding results of Schupp were
eventually subsumed by the following result of
Ol{\cprime}shanski\u{\i} (this result also subsumes our Theorem~A,
but not Theorem~B).

\begin{thmO}[Ol{\cprime}shanski\u{\i}~\cite{olshanskii:efficient}]
Let $|I| \geq 2$ and $\class= \{ H_i \mid i \in I\}$ be a countable
family of countable nontrivial groups.

There exists a 2-generated simple group $S$ that contains an
isomorphic copy of each member of $\class$ and, moreover, has the
following properties (in what follows, the copy of $H_i$ in $S$ is
denoted by $H_i$).

\textup{(1)} If $i,j \in I$, $i \neq j$, $|H_j| \geq 3$, then $S$ is
generated by $H_i$ and $H_j$.

\textup{(2)} If $i,j \in I$, $i \neq j$, then $H_i \cap H_j = 1$.

\textup{(3)} Every element of finite order in $S$ is conjugate to an
element in $H_i$, for some $i \in I$.

\textup{(4)} Every proper subgroup of $S$ is either infinite cyclic,
or infinite dihedral, or it conjugate of a subgroup of $H_i$, for
some $i \in I$.

\textup{(5)} If, for some $i \in I$, $x \in H_i$, $x \neq 1$, $y
\not\in H_i$, then either $S$ is generated by $\{x,y\}$ or both $x$
and $y$ are involutions, or both $x$ and $xy$ are involutions,

\textup{(6)} If $i,j \in I$, $i \neq j$, then $H_i \cap H_j^x =
{1}$, for every element $x$ in $S$

\textup{(7)} For every $i \in I$, $H_i$ is malnormal in  $S$ (for
every $x \in S \setminus H_i$, $H_i \cap H_i^x = 1$).

\end{thmO}

Thus there is a natural trade off in our approach. We extend
Theorem~S of Schupp (by adapting his approach using small
cancelation theory) to arbitrary families of groups in a way that,
in the countable case, partially recovers Theorem~O of
Ol'shanski\u{\i}. A modest gain is achieved by the fact that the
taken approach allows us to handle families of groups that are not
necessarily countable. On the other hand, in the countable case, we
recover only a small subset of the conclusions that are obtained by
the more powerful (but also more onerous) graded diagram methods
introduced by Olshanski\u{\i}.

\section{Proofs and additional comments}

\begin{proof}[Proof of Theorem~A]
Reindex the family $\class$ (if necessary) so that it is indexed by
an initial segment $I$ of the set of natural numbers
$\N=\{0,1,2\dots\}$ (including the possibility $I=\N$, if $I$ is
infinite). Moreover, in case the cyclic group $C_2$ of order 2 is a
member of $\class$ set $H_0=C_2$ and make sure that this is the only
copy of $C_2$ in $\class$.

For each $i \in I$, embed $H_i$ into a 2-generated simple group
$S_i=\langle s_i,t_i \rangle$ (this can be done by Theorem~S) and
consider the free product $F=A*B*(*_{i \in I} S_i)$, where
$A=\langle a \mid a^2 \rangle = C_2$, $B=\langle b \mid b^3 \rangle
= C_3$.

For each index $i \in I$ define the words
\begin{align*}
 u_i &= (ab)^{(2i+1)n+n}\ (ab^{-1})\ (ab)^{(2i+1)n+n-1}\ (ab^{-1})\
 \dots\ (ab^{-1})\ (ab)^{(2i+1)n+1} \ s_i, \\
 v_i &= (ab)^{(2i+2)n+n}\ (ab^{-1})\ (ab)^{(2i+2)n+n-1}\ (ab^{-1})\
 \dots\ (ab^{-1})\ (ab)^{(2i+2)n+1} \ t_i,
\end{align*}
where $n$ is a positive integer to be specified at a later stage.

Choose a nontrivial element $h_0$ in $H_0$ and, for each $i>0$,
choose a pair of distinct nontrivial elements $h_i$ and $\bar{h}_i$
in $H_i$. For each pair of indices $i,j \in I$ with $0 \leq i< j $,
define the words
\begin{align*}
 w_{(a,i,j)} &= (h_ih_j)^n\ (h_i\bar{h}_j)\ (h_ih_j)^{n-1} \ (h_i\bar{h}_j)\
 \dots \ (h_i\bar{h}_j)\ (h_ih_j)^1 \ a \\
 w_{(b,i,j)} &= (h_ih_j)^{2n}\ (h_i\bar{h}_j)\ (h_ih_j)^{2n-1}\ (h_i\bar{h}_j)\
 \dots \ (h_i\bar{h}_j)\ (h_ih_j)^{n+1} \ b .
\end{align*}

Let $R$ be the set of words obtained by symmetrization (closure
under inversion and conjugation; see Remark~1 for a precise
definition) of the set of words
\[
 R' =
  \{\ w_{(a,i,j)}, \ w_{(b,i,j)} \mid i,j \in I, \ 0 \leq i < j \ \} \cup
  \{\ u_i, \ v_i, \ (h_ia)^n, \ (h_ib)^n \mid i \in I \ \}
\]
and let $H=\langle \ F \mid R \ \rangle$.

Choose $n$ that is relatively prime to 6 and is sufficiently large
to ensure that the set of words $R$ satisfies the small cancelation
condition $C'(1/6)$ over the free product $F=A*B*(*_{i \in I} S_i)$
(see Remark~1 for a definition of the small cancelation condition
over free products). It follows, by a result of
Lyndon~\cite[Theorem~IV]{lyndon:dehn-alg}
(see~\cite[Section~V.9]{lyndon-s:b-cgt} for an exposition), that all
factors in the free product $F$ are embedded in $H=\langle \ F \mid
R \ \rangle$.

The $u$ relators and the $v$ relators ensure that $H$ is generated
by $a$ and $b$. On the other hand, the $w$ relators ensure that $H$
is generated by $H_i$ and $H_j$ for any $i,j \in I$ with $0 \leq i
<j$.

Let $M$ be a maximal normal subgroup of $H$ and let $S=H/M$. The
group $S$ is simple by the maximality of $M$. We claim that all the
factors $S_i$, $ i \in I$, are still embedded in $S$. The factor
$S_i$, being simple, either intersects $M$ trivially or is contained
in $M$. In the former case, the factor $S_i$ is still embedded in
$S=H/M$. The latter case implies that $h_i=1$ in $S$. Because of the
relators $(h_ia)^n$ and $(h_ib)^n$, it follows that $a^n=b^n=1$ in
$S$. However, $n$ is chosen to be relatively prime to 6. Thus
$a=b=1$ in $S$, which means that $S$ is trivial, a contradiction.

This completes the proof.
\end{proof}

We note here the crucial role of the embeddings $H_i \hookrightarrow
S_i$ in the course of the proof. On one hand, the number of
generators needed for each factor in $*_{i \in I} S_i$ is
uniformized. This is notationally convenient, but not crucial. More
significant is the simplicity of the factors $S_i$, which, helped by
the relators $(h_ia)^n$ and $(h_ib)^n$, ``protects'' the embedded
subgroups $H_i$ from ``crashing'' when $M$ is factored out from $H$.

\begin{proof}[Proof of Corollary~A]
Apply Theorem~A to $\class = \{ H_i \mid i \geq 1 \}$, where
$H_0=C_2$, $H_1=G$, and $H_{2i-4} = H_{2i-3}= C_i$, for $i \geq 3$
($C_m$ denotes the cyclic group of order $m$).
\end{proof}

\begin{proof}[Proof of Theorem~B]
Let $J$ be an indexing set of cardinality $\lambda$. For each $i \in
I$, embed $H_i$ into a simple group $S_i=\langle \{ s_{i,j} \mid j
\in J \}\rangle$. The cardinality of the simple group $S$ and the
generating system of $S_i$ can be chosen to be equal to
$\lambda=|J|$ by Theorem~S. Consider the free product $F=A*B*(*_{i
\in I} S_i)$, where $A=\langle \ a \mid a^2 \ \rangle =C_2$, $B=
*_{j \in J} \langle \ b_j \mid b_j^3 \ \rangle= *_{j \in J} C_3$.

Let $\alpha: I \times J \to J$ be an injective map (such a map
exists since $|I|\leq |J|$ and $|J|$ is infinite).

For each pair $(i,j) \in I \times J$, define the word
\[
 u_{i,j} = (ab_{\alpha(i,j)})^n \ (ab_{\alpha(i,j)}^{-1})\
 (ab_{\alpha(i,j)})^{n-1} \ (ab_{\alpha(i,j)}^{-1})\ \dots \
 (ab_{\alpha(i,j)}^{-1})\ (ab_{\alpha(i,j)})^1 \ s_{i,j} , \\
\]
where $n$ is a positive integer to be specified at a later stage.

For each $i \in I$, choose a nontrivial element $h_i$ in $H_i$. For
each $i' \in I$ such that $|H_{i'}|=\lambda$, choose a nontrivial
element $\bar{h}_{i'}$ in $H_{i'}$ different from $h_{i'}$ and
distinct nontrivial elements $h_{i',j}$, $j \in J$,
$\bar{h}_{i',j}$, $j \in J$, in $H_{i'}$ that are also different
from $h_{i'}$ and $\bar{h}_{i'}$. Let $L \subseteq I \times I$ be a
set of pairs such that, for each pair of indices $i,i' \in I$ such
that $|H_i|<|H_{i'}|=\lambda$, the ordered pair $(i,i')$ is in $L$,
and, for each pair of indices $i,i' \in I$ such that $i \neq i'$ and
$|H_i|=|H_{i'}|=\lambda$, exactly one of the ordered pairs $(i,i')$
and $(i',i)$ is in $L$ (if more than one member of $\class$ has
cardinality $\lambda$, then there are many possible choices for $L$
and we select one; the set $L$ must be nonempty because $|I| \geq
2$). For every pair $(i,i')$ in $L$, define the words
\begin{align*}
 w_{(a,i,i')} &= (h_ih_{i'})^n\ (h_i\bar{h}_{i'})\ (h_ih_{i'})^{n-1}\
  (h_i\bar{h}_{i'})\ \dots \ (h_i\bar{h}_{i'})\ (h_ih_{i'})^1 \ a \\
 w_{(b_j,i,i')} &= (h_ih_{i',j})^n\ (h_i\bar{h}_{i',j})\
   (h_ih_{i',j})^{n-1}\ (h_i\bar{h}_{i',j})\ \dots \
   (h_i\bar{h}_{i',j})\ (h_ih_{i',j})^1 \ b_j, \ j \in J.
\end{align*}

For all $i \in I$ and $j \in J$ also define the words
\[ (h_ia)^n, \ (h_ib_j)^n. \]

The group $H=F/R$ is defined as before ($R$ is obtained by
symmetrization of the set of all words defined so far and $n$ is
chosen to be relatively prime to 6 and to be sufficiently large to
yield the cancelation condition $C'(1/6)$ over the free product
$F=A*B*(*_{i \in I} S_i)$). The group $S=H/M$, where $M$ is a
maximal normal subgroup of $H$ satisfies the required conditions.
\end{proof}

\begin{remark}
We recall here the definition of the small cancelation property
$C'(1/6)$ over free products and specify a value for $n$ in the
proof of Theorem~A that ensures that this condition is satisfied.

Let $G=*_{i \in I} G_i$ be a free product of a nonempty family of
nontrivial groups $G_i$, $i \in I$ (here $I$ is an arbitrary
nonempty indexing set). The free product $G$ is generated by the set
$\Sigma=\cup_{i \in I} G_i \setminus \{1\}$ of nontrivial elements
in the (disjoint) union of the factors of $G$. A word $g_1g_2 \dots
g_k$ over $\Sigma$ is reduced if, for $\ell=1,\dots,k-1$, $g_\ell$
and $g_{\ell+1}$ come from a different factor of $G$. Every element
$g$ in $G$ can be represented by a unique reduced word over
$\Sigma$, called the normal form of $g$. By definition, the length
of an element in $G$ is equal to the length of its normal form. A
reduced word $g_1g_2 \dots g_k$ over $\Sigma$ is weakly cyclically
reduced if $k \leq 1$ or $g_kg_1 \neq 1$ (thus $g_k$ and $g_1$ may
come from the same factor of $G$, but may not cancel). Let
$u=u_1\dots u_k$ and $v=v_1\dots v_m$ be two reduced words over
$\Sigma$. If $k=0$ or $m=0$ or $u_kv_1 \neq 1$, we say that the
product $uv$ is semi-reduced (thus $u_k$ and $v_1$ may come from the
same factor of $G$, but may not cancel). A set $R$ of words over
$\Sigma$ is called symmetrized if it consists of weakly cyclically
reduced words, it is closed for inversion, and, for every $r$ in
$R$, all weakly cyclically reduced words over $\Sigma$ representing
conjugates of $r$ are in $R$. A nonempty reduced word $p$ over
$\Sigma$ is a piece in the symmetrized set $R$ if there exist two
reduced words $q_1$ and $q_2$ and two distinct words $r_1$ and $r_2$
in $R$ such that $r_1=pq_1$, $r_2=pq_2$ in $G$ and the products
$pq_1$ and $pq_2$ are semi-reduced. In this case we say that $p$ is
a piece in $r_1$ (and in $r_2$). Note that $p$ does not have to be a
subword of $r_1$ to be a piece in it. A symmetrized set $R$ of words
over $\Sigma$ satisfies the small cancelation property $C'(1/6)$
over the free product $G$ if, every word in $R$ has length greater
than 6 and, for every piece $p$, every reduced word $q$, and every
word $r$ in $R$ such that $r=pq$ in $G$ and the product $pq$ is
semi-reduced, the inequality $|p|<\frac{1}{6}|r|$ holds.

We now go back to our concrete setup from the proof of Theorem~A.
The lengths of the words in the set of relators $R'$ are
\begin{alignat*}{5}
 |u_i| &= 4n^2i +3n^2+3n-1, \qquad &&|w_{(a,i,j)}| &&= n^2+3n-1,
 \qquad &&|(h_ia)^n|&&=2n,\\
 |v_i| &= 4n^2i +5n^2+3n-1, \qquad &&|w_{(b,i,j)}| &&= 3n^2+3n-1,
 \qquad &&|(h_ib)^n|&&=2n.
\end{alignat*}
The reduced word $p = (h_ih_j)^n \ (h_i\bar{h}_j) (h_ih_j)^n$ of
length $4n+2$ is a piece in $w_{(a,i,j)}$, since it is a subword of
$w_{(b,i,j)}$ (and whence a prefix in a cyclic conjugate of
$w_{(b,i,j)}$) and since $w_{(a,i,j)}$ can be written as a
semi-reduced product $pq$ (for an appropriately chosen reduced word
$q$ with first letter $h_j^{-1}\bar{h}_j$). It can be easily
verified that $w_{(a,i,j)}$ does not have a piece longer than $p$
(although it has other pieces of the same length). Thus $n$ needs to
be selected in such a way that
\[ \frac{4n+2}{n^2+3n-1} < \frac{1}{6}. \]
This is true for any $n \geq 22$, but since we require $n$ to be
relatively prime to $6$, the smallest good choice is $n=23$.

We may now fix $n=23$, consider all other pieces of words, and
verify that the $C'(1/6)$ condition is satisfied.

For instance,  the word $(ab)^{(2i+1)n+n}(ab^{-1})(ab)^{(2i+1)n+n}$
of length $8ni+8n+2$ is a piece of $u_i$ (and this word does not
have any longer pieces). Thus we need to verify that
\[ \frac{8ni + 8n + 2}{4n^2i +3n^2+3n-1} < \frac{1}{6}, \]
for all $i \geq 0$. Think of the fraction on the left as a function
of $i$. Since $8n(3n^2+3n-1) - (8n+2)\cdot4n^2 <0$, this function is
decreasing for $i \geq 0$, the maximum is achieved at $i=0$, and its
value is $(8n+2)/(3n^2+3n-1) = 186/1655 < 1/6$.

We can equally easily verify all other cases. Thus we may take
$n=23$.
\end{remark}

\begin{remark}
Consider again the proof of Theorem~A. We used the original work of
Schupp not only to model our approach, but also to embed each group
$H_i$ into a simple 2-generated group $S_i$ (in order to protect
$H_i$ in the quotient $S=H/M$). In turn, in his proof of Theorem~S,
Schupp uses embeddings of $G$, $H$, and $K$ into countable simple
groups. At about the same time Schup proved his result, Goryushkin
also proved that every countable group can be embedded into a
2-generated simple group~\cite{goryushkin:simple}. Before the
results of Schupp and Goryushkin, it was known from the work of
P.~Hall that every countable group can be embedded into a
3-generated simple group~\cite[Theorem~C2]{phall:join}. However,
both Hall and Goryushkin also base their proofs on the existence of
embeddings of countable groups into countable simple groups. Thus to
get back on some firm footing one could perhaps go back directly to
the classical embedding results of Higman, Neumann and Neuman.
Namely they prove~\cite{higman-n-n:hnn} that every countable group
can be embedded into a countable group in which any two elements
that have the same order are conjugate. As a corollary, every
countable group can be embedded into a countable simple divisible
group (see~\cite[Theorem~IV.3.4]{lyndon-s:b-cgt} for an exposition).
Of course, by using such embeddings directly in the course of the
proof of Theorem~A, we could skip over a layer in the construction
at the cost of a mild notational difficulty (one would have to deal
with countably many countable generating sets).
\end{remark}

\subsection*{Acknowledgments}
The author would like to thank Centre Interfacultaire Bernoulli at
EPF-Lausanne for the support, the staff members for the hospitality
during his stay in May 2007, Goulnara Arzhantseva and Alain Valette
for their kind invitation to participate in the program, and the
referee for his/her rather useful remarks.


\def\cprime{$'$}

\end{document}